\documentclass[12pt]{article}
\usepackage{amsmath, amssymb}
\usepackage{epsfig}

\begin{document}

\title{
Concerning Riemann Hypothesis}

\author{R. Acharya
\footnote{I dedicate this note to my teacher George Sudarshan and
to the memory of Srinivasa Ramanujan (22 December 1887 -- 26 April
1920), ``The Man Who Knew `Infinity'."} \\
Department of Physics \& Astronomy \\
Arizona State University\\
Tempe, AZ 85287-1504}

\date{March 17, 2009}

\maketitle

\begin{abstract}
We present a quantum mechanical model which establishes the
veracity of the Riemann hypothesis that the non-trivial zeros of
the Riemann zeta-function lie on the critical line of $\zeta(s)$.
\end{abstract}

\begin{center}
\fbox{
\scalebox{0.90}{\includegraphics{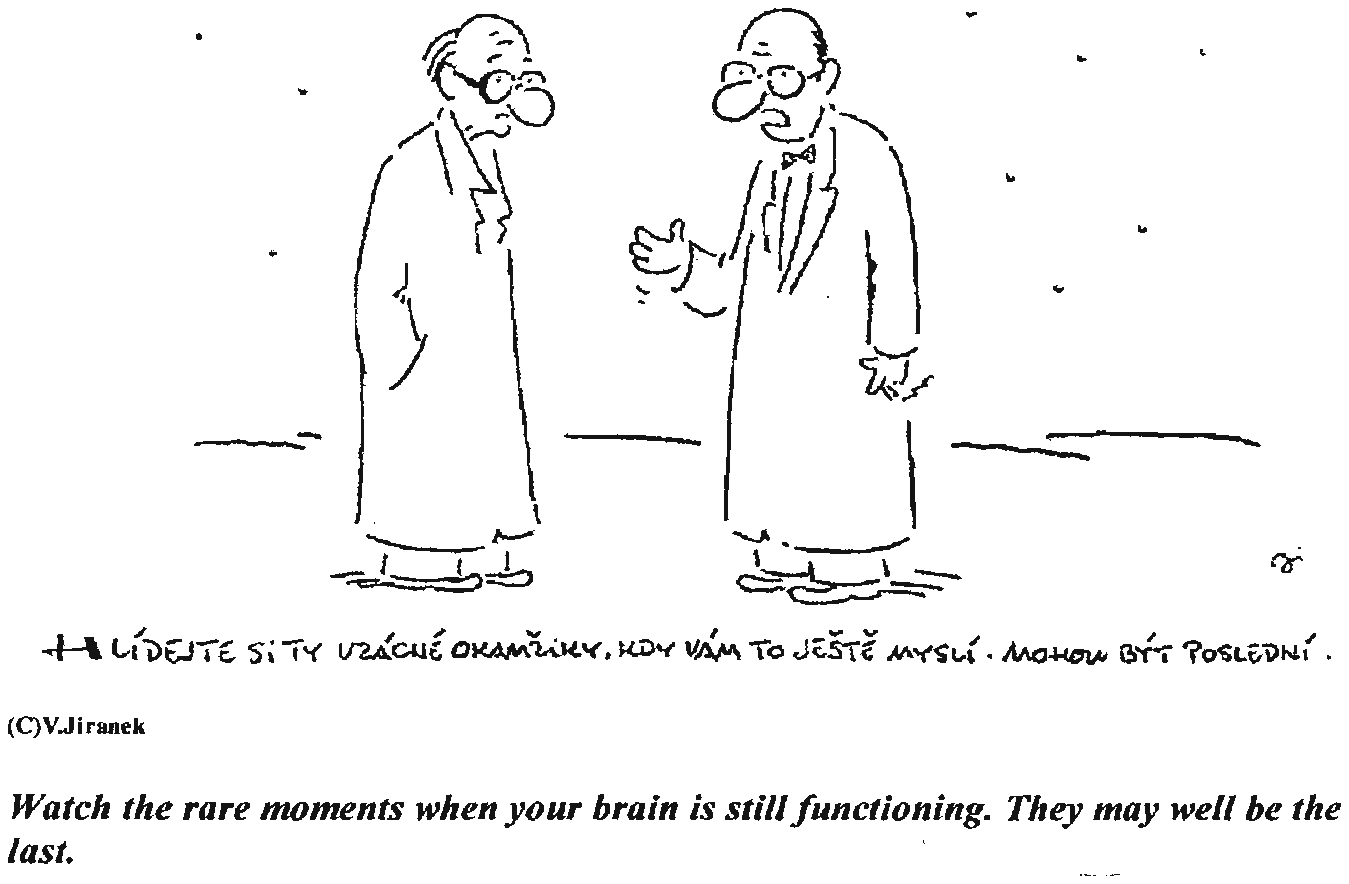}}
}
\end{center}

\vspace{0.25in}

We begin by recalling the all-too-familiar lore that the Riemann
hypothesis has been the Holy Grail of mathematics and physics for
more than a century \cite{Borwein}. It asserts that \textbf{all} the zeros of
$\xi(s)$ have $\sigma = \frac{1}{2}$, where $s=\sigma \pm it_{n}$,
$n = 1, 2, 3 \dots \infty$. It is believed all zeros of $\xi(s)$
are simple. The function $\zeta(s)$ is related to the Riemann
$\xi(s)$ function via the defining relation \cite{Borwein},

\begin{equation}\label{eq:1}
\xi(s) = \frac{1}{2} s(s-1)\pi^{-\frac{S}{2}} \Gamma \Bigl( \frac{s}{2} \Bigr)
\zeta(s)
\end{equation}
so that $\xi(s)$ is an entire function, where

\begin{equation}\label{eq:2}
\zeta(s) = \sum^{\infty}_{n=1} n^{-s}, \quad s = \sigma + it, \quad \sigma >
1
\end{equation}

$\zeta(s)$ is holomorphic for $\sigma > 1$ and can have No zeros for
$\sigma > 1$. Since $1/\Gamma(z)$ is entire, the function $\Gamma
\bigl( \frac{s}{2} \bigr)$ is non-vanishing, it is clear that
$\xi(s)$ \textbf{also} has no zeros in $\sigma > 1$: the zeros
of $\xi(s)$ are confined to the ``critical strip" $0 \leqslant \sigma
\leqslant 1$. Moreover, if $\rho$ is a zero of $\xi(s)$, then so is
$1-\rho$ and since $\overline{\xi(s)} = \xi(\overline{s})$, one deduces that
$\overline{\rho}$ and $1 - \overline{\rho}$ are also zeros. Thus the Riemann
zeros are symmetrically arranged about the real axis
\textbf{and} also about the ``critical line" given by $\sigma =
\frac{1}{2}$. The Riemann Hypothesis, then, asserts that ALL zeros
of $\xi(s)$ have Re $s= \sigma = \frac{1}{2}$.

We conclude this introductory, well-known remarks with the
assertion that every entire function $f(z)$ of \textbf{order
one} and ``\textbf{infinite type}" (which guarantees the
existence of \textbf{infinitely many Non-zero}
zeros can be represented by the Hadamard factorization, to wit \cite{Titchmarsh},

\begin{equation}\label{eq:3}
f(z) = z^{m} e^{A} e^{BZ} \prod^{\infty}_{n=1}
\Bigl( 1-\frac{z}{z_{n}} \Bigr) \exp \Bigl( \frac{z}{z_{n}} \Bigr)
\end{equation}
where `$m$' is the multiplicity of the zeros (so that $m=0$, for
simple zero).

Finally, $\xi(s) = \xi(1-s)$ is indeed an entire function of order
one and infinite type and it has \textbf{No zeros either for
$\mathbf{\sigma > 1}$ or $\mathbf{\sigma < 0}$.}

We invoke the well-known results of Lax-Phillips \cite{Lax} and
Faddeev-Pavlov \cite{Faddeev} scattering theory of automorphic functions in the
Poincare' upper-half plane, $z = x + iy$, $y > 0$, $ -\infty < x <
\infty$, which was motivated by Gelfand's \cite{Gelfand} observation of the
analogy between the Eisenstein functions \cite{Iwaniec} and the scattering
matrix, $S(\lambda)$. Recently, Yoichi Uetake \cite{Uetake} has undertaken a
detailed study of Lax-Phillips and Faddeev-Pavlov analysis, by
resorting to the technique of Eisenstein transform.

We recast this result by specializing to the case where incoming
and outgoing subspaces are \textbf{necessarily orthogonal}, $D_{-}
\perp D_{+}$ and by restricting to the case of Laplace-Beltrami
operator with \textbf{constant} `$x$' and identifying the
resulting wave equation, as a \textbf{one-dimensional} Schrodinger
equation, at \textbf{zero energy}, with a repulsive,
inverse-square potential, $V(y) = \lambda y^{-2}$, $\lambda > 0$,
$y > 0$ \textbf{and} an \textbf{infinite barrier at the origin},
$V(y) = \infty$, $y \leqslant 0$. We obtain zero-energy Jost
\cite{Jost,Chadan} function $F_{+} (k^{2} = 0)$:

\begin{equation}\label{eq:4}
S^{-1} (k^{2} = 0) = F_{+} (k^{2} = 0) = \frac{\xi(-2s)}{\xi(2s)}
\end{equation}
where $\xi(s)$ is Riemann's $\xi$ function \cite{Borwein}, Eq.(\ref{eq:1}) and we have shifted the variable `$s$' by
$\frac{1}{2}$, following the convention of Uetake \cite{Uetake}.
Since \textbf{all} zeros of Jost function $F_{+}(s)$ lie on the
critical line, $R s = - \frac{1}{4}$, we conclude that the Riemann
hypothesis is valid.

As an introduction to set the notation, we begin by presenting the
familiar case of the Euclidean plane,

\begin{equation}\label{eq:5}
\mathbb{R}^{2} = \{ (x,y) : x,y \in \mathbb{R} \}
\end{equation}

The group $G = \mathbb{R}^{2}$ acts on itself as translations, and
it makes $\mathbb{R}^{2}$ a homogeneous space. The Euclidean plane
is identified by the metric

\begin{equation}\label{eq:6}
d L^{2} = dx^{2} + dy^{2}
\end{equation}
with zero curvature $(K = 0 )$ and the Laplace-Beltrami operator
associated with this metric is given by

\begin{equation}\label{eq:7}
D = \frac{\partial^{2}}{\partial x^{2}} + \frac{\partial^{2}}{\partial y^{2}}
\end{equation}

Clearly, the exponential functions

\begin{equation}\label{eq:8}
\varphi (x,y) =
\exp[2\pi i (ux+vy)], \quad (u,v) \in \mathbb{R}^{2}
\end{equation}
are eigenfunctions of $D$:

\begin{equation}\label{eq:9}
(D + \lambda) \varphi = 0, \quad
\lambda = \lambda(\varphi) = 4\pi^{2} (u^{2} + v^{2})
\end{equation}

The upper-half hyperbolic half-plane \cite{Lax} (called the Poincare' plane)
is identified by

\begin{equation}\label{eq:10}
|H| = \{ z = x + iy \quad : \quad x \in \mathbb{R}, \quad
y \in \mathbb{R}^{+} \}
\end{equation}

$|H|$ is a Riemann manifold with the metric

\begin{equation}\label{eq:11}
d L^{2} = y^{-2} (dx^{2} + dy^{2})
\end{equation}

It represents a model of non-Euclidean geometry, where the role of
non-Euclidean motion is taken by the group $G$ of fractional
linear transformations,

\begin{equation}\label{eq:12}
W \rightarrow \frac{aw + b}{cw + d}
\end{equation}

\begin{equation}\label{eq:13}
\text{with} \quad ab - bc = 1, a, b,c, d \quad \text{real}
\end{equation}

The matrix $\left(
\begin{array}{cc}
  a & b \\
  c & d \\
\end{array}
\right)$
and its negative furnish the same transformation. The Riemannian
metric, Eq.(\ref{eq:11}) is invariant under this group of motions.

The Laplace-Beltrami operator is given by

\begin{equation}\label{eq:14}
y^{2} \Delta = ( \frac{\partial^{2}}{\partial x^{2}} +
\frac{\partial^{2}}{\partial y^{2}} )
\end{equation}

A discrete subgroup of interest is the modular group consisting of
transformations with integer $a$, $b$, $c$, $d$.

A fundamental domain $\mathbb{F}$ for a discrete subgroup $\Gamma$
is a subdomain of the Poincare' plane such that every point of $\Pi$ can
be carried into a point of the closure $\overline{\mathbb{F}}$ of
$\mathbb{F}$ by a transformation in $\Gamma$ and no point of
$\mathbb{F}$ is carried into another point of $\mathbb{F}$ by such
a transformation. $\overline{\mathbb{F}}$ can be regarded as a
manifold where those boundary points which can be mapped into each
other by a $\gamma$ in $\Gamma$ are identified.

Then, a function $f$ defined on $\Pi$ is called automorphic with
respect to a discrete subgroup $\Gamma$ if

\begin{equation}\label{eq:15}
f(\gamma w) = f(w)
\end{equation}
for all $\gamma$ in $\Gamma$.

By virtue of Eq.(\ref{eq:15}), an automorphic function is
completely determined by its values on $\overline{\mathbb{F}}$.
The Laplace-Beltrami operator, Eq.(\ref{eq:14}) maps automorphic
functions into automorphic functions.

In regular coordinates, $z = z + iy$, if we require $f(z)$ to be a
function of $y$ only, i.e., constant in $x$, we arrive at

\begin{equation}\label{eq:16}
\frac{\partial^{2}f(y)}{\partial y^{2}} +
\frac{\lambda_{0}f(y)}{y^{2}} = 0, \quad y > 0
\end{equation}
with two independent solutions \cite{Iwaniec},

\begin{equation}\label{eq:17}
\frac{1}{2} ( y^{s} + y^{1-s} ) \quad \text{and} \quad
\frac{1}{2s-1} ( y^{s} - y^{1-s} )
\end{equation}
where
\begin{equation}\label{eq:18}
\lambda_{0} = s(1-s).
\end{equation}

For $ s = \frac{1}{2}$ ($\lambda_{0} = \frac{1}{4}$), the above
eigenfunctions become $y^{\frac{1}{2}}$ and $y^{\frac{1}{2}}
\log y$ respectively.

We now make an \textbf{important} observation which is \textbf{crucial},
i.e., we can view Eq.(\ref{eq:16}) as a Schrodinger equation at
zero energy, for an inverse-square potential, i.e.,

\begin{equation}\label{eq:19}
\frac{\partial^{2}\Psi(y)}{\partial y^{2}} + (k^{2} - V(y))
\Psi(y) = 0
\end{equation}
where
\begin{equation}\label{eq:20}
k^{2} = 0 \left( \frac{2m}{\hbar^{2}} = 1 \right) \quad \text{and} \quad V(y) =
\frac{s(s-1)}{y^{2}}
\end{equation}
with the \textbf{all-important constraint} that
\begin{gather}
V(y) = \frac{s(s-1)}{y^{2}}, \quad y > 0  \label{eq:21}\\
[ s > 1 \quad \text{or} \quad s < 0] \quad \text{i.e. Repulsive!} \notag
\end{gather}
and
\begin{equation}\label{eq:22}
V(y) = \infty, \quad y \leqslant 0.
\end{equation}

In other words, the restriction on variable `$y$' in the Poincare'
upper-half plane, Eq.(\ref{eq:10}) requires that $ y > 0 $ and
Eq.(\ref{eq:22}) ensures that this requirement is satisfied by
placing an ``infinite barrier" at $y=0$, so that $\Psi(y) \equiv
0$, in the ``left-half" $y$-axis.

We now follow Uetake's analysis closely \cite{Uetake}.The
Eisenstein series of two variables $E(z,s)$ on a fundamental
domain is by definition, is real analytic for $z=x+iy \in |H|$ and
$E(z,s)$ is meromorphic in $s$ in the complex plane $C$. It is
regular for each $z \in |H|$, with respect to `$s$' in $R s
\geqslant 0$, except at $s = \frac{1}{2}$. In fact, $E(z,
\frac{1}{2} + s)$ is an automorphic function on $|H|$:

\begin{equation}\label{eq:23}
E \left( \gamma(z), \frac{1}{2} + s \right) = E \left( z, \frac{1}{2} + s
\right) \quad \forall \quad \gamma \in SL_{2}(z).
\end{equation}

In Eq.(\ref{eq:23}), following Uetake, we have shifted the
variable `$s$' by $\frac{1}{2}$ and $\gamma(z) = z + n$ for some
$n \in Z$. Thus, $E(z,\frac{1}{2} + s)$ can be viewed as a
function on $|H|$ (upper-half Poincare' plane.)

Furthemore, the Eisenstein series is a (non $L^{2} - $)
eigenfunction of the non-Euclidean Laplacian, $-y^{2} \left(
\frac{\partial^{2}}{\partial x^{2}} + \frac{\partial^{2}}{\partial
y^{2}} \right)$ on $|H|$, i.e.,

\begin{equation}\label{eq:24}
\left[ -y^{2} \left( \frac{\partial^{2}}{\partial x^{2}} + \frac{\partial^{2}}{\partial
y^{2}} \right) - \frac{1}{4} \right] E(z, \frac{1}{2} + i t) =
t^{2} E(z, \frac{1}{2} + i t )
\end{equation}
for all $t \in C$. Integrating Eq.(\ref{eq:24}), w.r.t x yields Eq.(\ref{eq:19}) and
Eq.(\ref{eq:20}) where $s \rightarrow s + \frac{1}{2}$ \cite{Lax}.

Eq.(\ref{eq:24}) is an eigenfunction property of the Eisenstein
series.

Also, the Eisenstein series satisfies the functional equation,

\begin{equation}\label{eq:25}
E \left( z, \frac{1}{2} + s \right) = S(s) E  \left( z,
\frac{1}{2}- s \right)
\end{equation}
where $S(s)$ is the scattering ``matrix." Eq.(\ref{eq:24}) and
Eq.(\ref{eq:25}) are proven, for instance, in Y. Motohashi's text
on the spectral theory of the Riemann Zeta-function
\cite{Motohashi}.

When $E \left( z, \frac{1}{2} + s \right) $ is expanded in the
Fourier series of $\exp(inz)$, the zero Fourier coefficient takes
the form,

\begin{equation}\label{eq:26}
y^{\frac{1}{2} + s} + S(s) y^{\frac{1}{2} - s}
\end{equation}
where the scattering matrix $S(s)$ has the following form, for
\textbf{orthogonal} incoming and outgoing subspaces $D_{-}$ and
$D_{+}$ \cite{Uetake}:

\begin{equation}\label{eq:27}
S(s) = \frac{\xi (2s)}{\xi (-2s)}
\end{equation}
where Riemann's $\xi(s)$ function is related to THE Riemann
function $\zeta(s)$ via Eq.(\ref{eq:1}) \cite{Titchmarsh}:

\begin{equation*}\label{eq:1a}
\xi(s) = \frac{1}{2} \pi^{-\frac{S}{2}} s(s-1) \Gamma \left( \frac{s}{2} \right)
\zeta(s) \tag{1}
\end{equation*}

Importantly, $S(s)$ has poles in $ -\frac{1}{2} < R s < 0 $ and the
Riemann hypothesis is the assertion that \textbf{all} poles of
$S(s)$ lie on $ Rs = - \frac{1}{4}$.

We move on next to make contact with the Jost functions
\cite{Jost} $F_{\stackrel{+}{-}}(k)$ and the $S$ matrix.
It is well-known that the Jost function is identified via the
defining relation,

\begin{equation}\label{eq:28}
S(k) = \frac{F_{-}(k)}{F_{+}(k)}
\end{equation}
where $F_{+}(k)$ is holomorphic in the upper-half complex $k$-plane
and $F_{-}(k)$ is holomorphic in the lower-half $k$-plane, where
$k$ is the momentum.

Thus,
\begin{equation}\label{eq:29}
F_{+}(k) = S^{-1}(k) F_{-}(k), \quad Imk = 0, \quad -\infty < k <
\infty
\end{equation}
and

\begin{equation}\label{eq:30}
\det S(k) \neq 0, \quad Imk = 0
\end{equation}

The ``solution" to the boundary value problem was formulated by
Krutov, Muravyev and Troitsky in 1996 \cite{Krutov}. It reads:

\begin{equation}\label{eq:31}
F_{\stackrel{+}{-}}(k) = \Pi_{\stackrel{+}{-}}(k) \exp
\left( \frac{1}{2 \pi i} \int^{\infty}_{-\infty}
\frac{\ln(S^{-1}(k) \pi^{2}_{-}(k)}{k' - k \mp io}dk'
\right)
\end{equation}
where
\begin{equation}\label{eq:32}
\Pi_{\stackrel{+}{-}}(k) = \prod^{m'}_{j=1} \frac{k\mp ik_{j}}{k\pm
ik_{j}}, \quad k_{j} > 0, \quad j=1,2\ldots m' (m' < \infty)
\end{equation}
and
\begin{equation}\label{eq:33}
\Pi_{\stackrel{+}{-}}(k) \equiv 1, \quad m' = 0
\end{equation}
where $m'$ is the number of bound states.

It is, of course, well known that

\begin{equation}\label{eq:34}
\frac{1}{k'- k \mp io} = P \frac{1}{k' - k} \pm i \pi \delta
(k'-k)
\end{equation}

Proceeding further, \textbf{we set} $\mathbf{k=0}$ \textbf{(zero
energy).} The simplification is immediate. One finds that the zero
energy Jost functions are given by:

\begin{equation}\label{eq:35}
S^{-1}(s) = F_{+}(s), \quad F_{-}(s) = 1
\end{equation}

Thus, we arrive at the stated result, Eq.(\ref{eq:14}):

\begin{equation*}\label{eq:14a}
S^{-1} (s) = F_{+} (s) = \frac{\xi(-2s)}{\xi(2s)}
\tag{4}
\end{equation*}

The Riemann hypothesis is the assertion that \textbf{all poles} of
$S(s)$ lie on $Rs = - \frac{1}{4}$. From Eq.(\ref{eq:4}), we see
that this requires that \textbf{all zeros} of the Jost function,
$F_{+}(s)$ must lie on $Rs = - \frac{1}{4}$.

It is well-known from the the theory of Jost functions, summarized
by Khuri \cite{Khuri} that the $s$-wave Jost function for a
potential is \textbf{an entire function of the ``coupling constant"}
$\mathbf{\lambda}$\textbf{, with an infinite number of zeros extending to
infinity.} For a repulsive potential $V(y)$ and at zero energy,
these zeros of $\lambda_{n}$ will all be real and negative,
$\lambda_{n}(0) < 0$. By changing variables to `$s$' (\textbf{in Uetake's
convention}), with $\lambda = s (s-1)$, it follows that as a
function of `$s$', the Jost function $F_{+}(s)$ has \textbf{only zeros} on
the line $S_{n} = - \frac{1}{4} \pm it$. (Again, in Uetake's
\cite{Uetake} convention for `$s$'!)

This is automatic, in view of Eq.({\ref{eq:21})! At this point, the
discussion in Khuri detailed in his Eq.(1.1), Eq.(1.2) and
Eq.(1.3) demonstrates the justification (required). While Khuri's
analysis \cite{Khuri} dealt with exponentially decreasing
potentials for $x \rightarrow \infty$, the present situation of
inverse-square potential can be dealt with by taking a slight
detour in arriving at the analogue of Khuri's Eq.(1.2), i.e., it
is easily demonstrated that $[k = i \tau]$:

\begin{equation}\label{eq:36}
I m \left[ \lambda_{n}(i \tau) \right] \int^{\infty}_{0} \left|
f(\lambda_{n} (i \tau); \text{ } i \tau, \text{ } y \right|^{2} dy = 0
\end{equation}

To see this, one writes down the differential equation for the
Jost \textbf{solution}, $f(\lambda; k; y)$ [which is identical to
the $s$-wave Jost function in 3 dimensions!] and its complex
conjugate partner, then, one can perform the ``canonical"
operation of multiplying the equation for $f(y)$ by $f^{*}(y)$ and
like-wise for $f^{*}(y)$ equation by $f(y)$, doing the subtraction
\textbf{and multiply through} by $V(y)$ {which is real!] Finally,
expressing the terms involving $ff^{*''} - f^{''}f^{*}$ as
derivative of the \textbf{wronskian} which is \textbf{independent of}
$y$ (!), one finally ends up with Eq.(\ref{eq:36}). It is an
elementary exercise to perform the integral in Eq.(\ref{eq:36}) by
noting that the Jost solution has the form:

\begin{equation}\label{eq:37}
f(k,y) = \sqrt{\frac{\pi k y}{2}} e^{i(\frac{\pi}{2} \nu +
\frac{\pi}{4})} H^{(1)}_{\nu} (ky)
\end{equation}
where $H^{(1)}_{\nu} (ky)$ is Hankel function of $I$ kind and it
has the required asymptotic behavior, i.e.,

\begin{equation}\label{eq:38}
f(k,y) \underset{y \rightarrow \infty}{\longrightarrow} e^{iky}
\end{equation}

One notes that \cite{Gradshteyn}

\begin{equation}\label{eq:39}
\int^{\infty}_{0} y K^{2}_{\nu(i \tau)} (y) dy = \frac{1}{8}
\frac{\pi \nu (i \tau)}{\sin \pi \nu (i \tau)} \text{, } \nu (i \tau) = \sqrt{ \lambda (i \tau) + \frac{1}{4}}
\end{equation}
where $\nu (i \tau) \neq 1, 2, 3, \ldots \infty$.

\begin{equation}\label{eq:40}
\left[ I m  \lambda_{n} (i \tau) \right] \frac{1}{\tau} \frac{2}{\pi} \int^{\infty}_{0} y K^{2}_{\nu (i \tau)} (y) dy = 0
\end{equation}

Since the integral in Eq.(\ref{eq:40}) is finite and non-vanishing
[keeping $\nu \neq 1, 2, 3\ldots \infty$], one ends up with (following Khuri), taking the limit $\tau \rightarrow 0$ and
get

\begin{equation}\label{eq:41}
I m \lambda_{n} (0) = 0
\end{equation}

The zero energy, coupling constant spectrum, $\lambda_{n}(0)$ \textbf{must}
lie on the negative real axis for $ V > 0$ (Repulsive potential).
The rationale behind the idea of ``getting rid of" the potential,
$V(y) = \frac{\lambda}{y^{2}}$ in Eq.(\ref{eq:36}) is to ensure
that the resulting equation (Eq.(\ref{eq:36})) is now \textbf{finite} and
one can then continue on to conclude that the zero energy,
coupling constant $\lambda_{n}(0)$ must lie on the negative real
line, for $V > 0$.

In conclusion, the derivation of Jost function, $F_{+} (s) =
\frac{\xi(-2s)}{\xi(2s)}$ where $\xi(s)$ is Riemann's $\xi$
function and the established key assertion that ALL zeros of
$F_{+} (s)$ must lie on the critical line, $Rs = -\frac{1}{4}$,
leads us to the conclusion that the Riemann's hypothesis is valid
\cite{Berry}.

\vspace{24pt}
\noindent \textbf{Acknowledgement:} I am indebted to Irina Long for her kind assistance.

\renewcommand\refname{References:}

\end{document}